\documentstyle[12pt]{article}
\newcommand{\calx}{\mbox{${\cal X}$}}
\newcommand{\calz}{\mbox{${\cal Z}$}}

\newcommand{\calo}{\mbox{${\cal O}$}}
\newcommand{\calp}{\mbox{${\cal P}$}}
\newcommand{\calk}{\mbox{${\cal K}$}}
\newcommand{\rama}{\mbox{ram$_{\cal X}(\alpha)$}}
\newcommand{\pa}{\mbox{$\partial_{C_i}(\alpha)$}}

\title{The $u$-invariant of the function fields of $p$-adic curves}

\author{R. Parimala and V. Suresh}

\date{}

\begin{document}

\maketitle

\begin{abstract}
The $u$-invariant of a field is the maximum dimension of
ansiotropic quadratic forms over the field. It is an open question
whether the $u$-invariant of function fields of $p$-aidc curves is
8. In this paper, we answer this question in the affirmative for
function fields of non-dyadic $p$-adic curves.
\end{abstract}

\section*{Introduction}

It is an open question ([L], Q. 6.7, Chap XIII)  whether every
quadratic form in at least nine variables over the function fields
of  $p$-adic curves has a non-trivial zero.  Equivalently, one may
ask whether the $u$-invariant of such a field is 8.  The
$u$-invariant of a field $F$ is defined as the maximal dimension
of anisotropic quadratic forms over $F$.  In this paper we answer
this question in the affirmative if the $p$-adic field is
non-dyadic.

In ([PS], 4.5), we  showed  that every quadratic form in eleven
variables over the function field of a $p$-adic curve, $p \neq 2$,
has a nontrivial zero. The main ingredients in the proof were the
following: Let $K$ be the  function field of a $p$-adic curve $X$
and $p \neq 2$.
\begin{enumerate}
\item{}(Saltman ([S1], 3.4)) Every element in the Galois
cohomology group $H^2(K, {\mathbf{Z}}/2{\mathbf{Z}})$ is a sum of
at most two symbols.

\item{} (Kato [K], 5.2) The unramified cohomology group $H^3_{\rm
nr}(K/{\cal X}, {\mathbf{Z}}/2{\mathbf{Z}}(2))$ is zero for a
regular projective model ${\cal X}$ of $K$.

\end{enumerate}

If $K$ is as above,  we proved ([PS],  3.9) that every element in
$H^3(K, {\mathbf{Z}}/2{\mathbf{Z}})$ is a symbol of the form
$(f)\cdot(g)\cdot(h)$ for some $f, g, h \in K^*$ and $f$ may be
chosen to be a value of a given binary form $<\! a, b \!>$ over
$K$. If further given  $\zeta = (f)\cdot(g)\cdot(h) \in H^3(K,
{\mathbf{Z}}/2{\mathbf{Z}})$ and a ternary form $<\! c, d, e \!>$,
one can choose $g', h' \in K^*$ such that $\zeta =
(f)\cdot(g')\cdot(h')$ with $g'$ a value of $<\! c,d,e \!>$, then,
one is led to the conclusion that $u(K)=8$ (cf 4.3). We in fact
prove that such a choice of $g',h' \in K^*$ is possible by proving
the following local global principle:

\paragraph*{Theorem.} Let $K=k(X)$ be the function field of a
curve $X$ over a $p$-adic field $k$. Let $l$ be a prime not equal
to p. Assume that $k$ contains a primitive $l^{\rm th}$ root of
unity. Given  $\zeta \in H^3(K, \mu_l^{\otimes 2})$ and  $\alpha
\in H^2(K, \mu_l)$ corresponding to a degree $l$ central division
algebra over $K$, satisfying $\zeta = \alpha \cup (h_v)$ for some
$(h_v) \in H^1(K_v, \mu_l)$, for all discrete valuations of $K$,
there exists $(h) \in H^1(K, \mu_l)$ such that  $\zeta = \alpha
\cup (h)$. In fact one can restrict the hypothesis to discrete
valuations of $K$ centered on codimension one points of a regular
model $\calx$, projective  over the ring of integers ${\cal O}_k$
of $k$.

A key ingredient towards the proof of the theorem is a recent
result of Saltman ([S3]) where the ramification pattern of prime
degree central simple algebras over function fields of $p$-adic
curves is completely described.

We thank J.-L. Colliot-Th\'el\`ene for helpful discussions during
the preparation of this paper and for his critical
comments on the text.

\section*{1. Some Preliminaries}

In this section we recall a few basic facts from the algebraic
theory of quadratic forms and Galois cohomology. We refer the
reader to ( [C]) and ([Sc]).

Let $F$ be a field and $l$ a prime not equal to the characteristic
of $F$. Let $\mu_l$ be the group of $l^{\rm th}$ roots of unity.
For $i \geq 1$, let $\mu_l^{\otimes i}$ be the Galois module given
by the tensor product of $i$ copies of $\mu_l$. For $n \geq 0$,
let $H^n(F, \mu_l^{\otimes i})$ be the $n^{\rm th}$ Galois
cohomology group with coefficients in $\mu_l^{\otimes i}$.

We have the Kummer isomorphism $F^*/F^{*^l} \simeq H^1(F, \mu_l)$.
For $a \in F^*$, its class in $H^1(F, \mu_l)$ is denoted by $(a)$.
If $a_1, \cdots , a_n \in F^*$, the cup product $(a_1) \cdots
(a_n) \in H^n(F, \mu_l^{\otimes n})$ is  called a {\it symbol}.
We have an isomorphism $H^2(F, \mu_l)$ with the $l$-torsion
subgroup $_lBr(F)$ of the Brauer group of $F$. We define the {\it
index} of an element $\alpha \in H^2(F, \mu_l)$ to be the index of
the corresponding central simple algebra in $_lBr(F)$.

Suppose $F$ contains all the $l^{\rm th}$ roots of unity. We fix a
generator $\rho$ for the cyclic group $\mu_l$ and identify the
Galois modules $\mu_l^{\otimes i}$ with ${\mathbf{Z}}/l{\mathbf{Z}}$. This
leads to an identification of $H^n(F, \mu_l^{\otimes m})$ with
$H^n(F, {\mathbf{Z}}/l{\mathbf{Z}})$. The  element in $H^n(F,
{\mathbf{Z}}/l{\mathbf{Z}})$ corresponding to the  symbol $(a_1)
\cdots (a_n) \in H^n(F, \mu_l^{\otimes n})$ through this
identification is again denoted by $(a_1) \cdots (a_n)$. In
particular  for $a, b \in F^*$, $(a) \cdot (b) \in H^2(K,
{\mathbf{Z}}/l{\mathbf{Z}})$ represents the cyclic algebra $(a,
b)$ defined by the relations $x^l = a$, $y^l = b$ and $xy = \rho
yx$.

Let $v$ be a discrete valuation of $F$. The residue field of $v$
is denoted by $\kappa(v)$. Suppose  char$(\kappa(v)) \neq l$. Then
there is a {\it residue} homomorphism $\partial_v : H^n(F,
\mu_l^{\otimes m}) \to H^{n-1}(\kappa(v), \mu_l^{\otimes (m-1)})$.
Let $\alpha \in H^n(F, \mu_l^{\otimes m})$. We say that $\alpha$
is {\it unramified} at $v$ if $\partial_v(\alpha) = 0$; otherwise
it is said to be {\it ramified} at $v$. If $F$ is complete
with respect to $v$, we denote the kernel of
$\partial_v$ by $H^n_{\rm nr}(F, \mu_l^{\otimes m})$.
Suppose $\alpha$ is unramified at $v$. Let $\pi \in
K^*$ be a parameter at $v$ and $\zeta = \alpha \cup (\pi) \in
H^{n+1}(F, \mu_l^{\otimes (m+1)})$. Let $\overline{\alpha} =
\partial_v(\zeta) \in H^n(\kappa(v), \mu_l^{\otimes m})$. The
element $\overline{\alpha}$ is independent of the choice of the
parameter $\pi$ and is called the {\it specialization} of $\alpha$
at $v$. We say that $\alpha$ {\it specializes to} $
\overline{\alpha}$ at $v$. The following result is well known.

\paragraph*{Lemma 1.1} Let $k$ be a field and $l$ a prime not
equal to the  characteristic of $k$. Let $K$ be a complete
discrete valuated field with residue field $k$. If  $H^n(k,
\mu_l^{\otimes 3}) = 0$ for $n \geq 3$,  then $H^3_{\rm nr}(K,
\mu_l^{\otimes 3}) = 0$. Suppose further that every element in
$H^2(k, \mu_l^{\otimes 2})$ is a symbol. Then every element in
$H^3(K, \mu_l^{\otimes 3})$ is a symbol.

\paragraph*{Proof.} Let $R$ be the ring of integers in $K$. The
Gysin exact sequence in \'etale cohomology yields an exact
sequence (cf. [C], p.21, \S 3.3)
$$
H^3_{\rm \acute{e}t}(R, \mu_l^{\otimes 3}) \to H^3(K,
\mu_l^{\otimes 3}) \buildrel{\partial}\over{\to} H^2(k,
\mu_l^{\otimes 2}) \to H^4_{\rm \acute{e}t}(R, \mu_l^{\otimes 3})
$$
Since $R$ is complete,  $ H^n_{\rm \acute{e}t}(R, \mu_l^{\otimes
3}) \simeq H^n(k, \mu_l^{\otimes 3})$ ([Mi], p.224 Corollary 2.7).
Hence $H^n_{\rm \acute{e}t}(R, \mu_l^{\otimes 3})= 0$ for $n \geq
3$, by the hypothesis. In particular $\partial : H^3(K,
\mu_l^{\otimes 3}) \to H^2(k, \mu_l^{\otimes 2})$ is an
isomorphism and $H^3_{\rm nr}(K, \mu_l^{\otimes 3}) = 0$. Let $u,
v \in R$ be units and $\pi \in R$ a parameter. Then we have
$\partial ((u) \cdot (v) \cdot (\pi)) = (\overline{u}) \cdot
(\overline{v})$. Let $\zeta \in H^3(K, \mu_l^{\otimes 3})$. Since
every element in $H^2(k, \mu_l^{\otimes 2})$ is a symbol, we have
$\partial(\zeta) = (\overline{u}) \cdot (\overline{v})$ for some
units  $u, v \in R$.  Since $\partial$ is an isomorphism, we have
$\zeta = (u) \cdot (v) \cdot (\pi)$. Thus every element in $H^3(K,
\mu_l^{\otimes 3})$ is a symbol. \hfill $\Box$

\paragraph*{Corollary 1.2} Let $k$ be a $p$-adic field and
$K$  the function field of an integral curve over $k$. Let $l$ be
a prime not equal to $p$. Let $K_v$ be the completion of $K$ at a
discrete valuation of $K$. Then $H^3_{nr}(K_v, \mu_l^{\otimes 3})
= 0$. Suppose further that $K$ contains a primitive $l^{\rm th}$
root of unity. Then every element in $H^3(K_v, \mu_l^{\otimes 3})$
is a symbol.

\paragraph*{Proof.} Let $v$ be a discrete valuation of $K$
and $K_v$ the completion of $K$ at $v$. The  residue field
$\kappa(v)$ at $v$ is either a $p$-adic field or a function field
of a curve over a finite field of characteristic $p$. In either
case, the cohomological dimension of $\kappa(v)$ is 2 and hence
$H^n(\kappa(v), \mu_l^{\otimes 3}) = 0$ for $n \geq 3$. By (1.1),
$H^3_{nr}(K_v, \mu_l^{\otimes 3}) = 0$.

If $\kappa(v)$ is a local field, by local class field theory,
every finite dimension central division algebra over $\kappa(v)$
is split by an unramified (cyclic) extension. If $\kappa(v)$ is a
function field of a curve over finite field, then by a classical
theorem of Hasse-Brauer-Neother-Albert, every finite dimensional
central division algebra over $\kappa(v)$ is split by a cyclic
extension. Since $\kappa(v)$ contains a primitive $l^{\rm th}$ of
unity, every element in $H^2(\kappa(v),
{\mathbf{Z}}/l{\mathbf{Z}})$ is a symbol.  By (1.1),  every
element in $H^3(K_v, {\mathbf{Z}}/l{\mathbf{Z}})$ is a symbol.
\hfill $\Box$

\vskip 5mm

Let $\calx$ be a regular integral scheme of dimension $d$, with
field of fractions $F$. Let $\calx^1$ be the set of points of
$\calx$ of codimension 1. A point  $x \in \calx^1$  gives rise to
a discrete valuation $v_x$ on $F$. The residue field of this
discrete valuation ring is denoted by $\kappa(x)$ or
$\kappa(v_x)$. The corresponding residue homomorphism is denoted
by $\partial_x$. We say that an element $\zeta \in H^n(F,
\mu_l^{\otimes m})$ is {\it unramified} at $x$ if
$\partial_x(\zeta) = 0$; otherwise it is said to be {\it ramified}
at $x$. We define the ramification divisor $\mbox{
ram}_{\calx}(\zeta) = \sum x$ as $x$ runs over $\calx^1$ where
$\zeta$ is ramified.  The unramified cohomology on $\calx$,
denoted by $H^n_{nr}(F/\calx, \mu_l^{\otimes m})$, is defined as
the intersection of kernels of the residue homomorphisms
$\partial_x : H^n(F, \mu_l^{\otimes m}) \to H^{n-1}(\kappa(x),
\mu_l^{\otimes (m-1)})$, $x$ running over $\calx^1$. We say that
$\zeta \in H^n(F, \mu_l^{\otimes m})$ is {\it unramified on
$\calx$} if $\zeta \in H^n_{nr}(F/\calx, \mu_l^{\otimes m})$. If
$\calx = {\rm Spec}(R)$, then we also say that $\zeta$ is
unramified on $R$ if it is unramified on $\calx$. Suppose $C$ is
an irreducible subscheme of ${\calx}$ of codimension 1. Then the
generic point $x$ of $C$ belongs to $\calx^1$ and we  set
$\partial_x$ = $\partial_C$. If $\alpha \in H^n(F, \mu_l^{\otimes
m})$ is unramified at $x$, then we say that $\alpha$ is {\it
unramified} at $C$.

Let $k$ be a $p$-adic field and $K$ the  function field of a
smooth, projective, geometrically integral curve $X$ over $k$. By
the resolution of singularities for surfaces (cf.  [Li1] and [Li2]),
there exists a regular, projective model $\calx$ of $X$ over the
ring of integers $\calo_k$ of $k$. We call such an $\calx$  a {\it
regular projective model of} $K$.  Since the generic fibre $X$ of
$\calx$ is  geometrically integral, it follows that the special
fibre ${\overline{\calx}}$ is connected. Further if $D$ is a
divisor on $\calx$, there exists a proper birational morphism
$\calx' \to \calx$ such that the total transform of $D$ on
$\calx'$ is a divisor with normal crossings (cf. [Sh], Theorem,
p.38 and Remark 2, p. 43). We use this result throughout this
paper without further reference.

Let $k$ be a $p$-adic field and $K$ the  function field of a
smooth, projective, geometrically integral curve over $k$. Let $l$
be a prime not equal to $p$. Assume that $k$ contains a primitive
$l^{\rm th}$ of unity. Let $\alpha \in H^2(K, \mu_l)$. Let $\calx$
be a regular  projective model of $K$ such that the ramification
locus ram$_{\calx}(\alpha)$ is a union of regular curves with
normal crossings.  Let $P$ be a closed point in the intersection
of two regular curves $C$ and $E$ in ram$_{\calx}(\alpha)$.
Suppose that $\partial_C(\alpha) \in
H^1(\kappa(C),{\mathbf{Z}}/l{\mathbf{Z}})$ and $\partial_E(\alpha)
\in H^1(\kappa(E), {\mathbf{Z}}/l{\mathbf{Z}})$ are unramified at
$P$. Let $u(P), v(P) \in H^1(\kappa(P),
{\mathbf{Z}}/l{\mathbf{Z}})$ be the specialisations at $P$ of
$\partial_C(\alpha)$ and $\partial_E(\alpha)$ respectively.
Following Saltman ([S3], \S 2), we say that $P$ is a {\it cool
point} if $u(P), v(P)$ are trivial and a {\it chilli point} if
$u(P)$ and $v(P)$ generate the same subgroup of $H^1(\kappa(P),
{\mathbf{Z}}/l{\mathbf{Z}})$ and neither of them is trivial. If
$u(P)$ and $v(P)$ do not generate the same subgroup of $
H^1(\kappa(P), {\mathbf{Z}}/l{\mathbf{Z}})$, then $P$ is said to
be a {\it hot point}. Let ${\cal O}_{\calx, P}$ be the regular
local ring at $P$ and  $\pi$, $\delta$ prime elements in ${\cal
O}_{\calx, P}$ which define $C$ and $E$ respectively at $P$. The
condition that $\partial_C(\alpha) \in H^1(\kappa(C),
{\mathbf{Z}}/l{\mathbf{Z}})$ and $\partial_E(\alpha) \in
H^1(\kappa(E), {\mathbf{Z}}/l{\mathbf{Z}})$ are unramified at $P$
is equivalent to the condition $\alpha = \alpha'+ (u, \pi) + (v,
\delta)$ for some units $u, v \in {\cal O}_{{\cal X}, P}$ and
$\alpha'$ unramified on ${\cal O}_{\calx, P}$  ([S2], \S 2). The
specialisations of $\partial_C(\alpha)$ and $\partial_E(\alpha)$
in $H^1(\kappa(P), {\mathbf{Z}}/l{\mathbf{Z}}) \simeq
\kappa(P)^*/\kappa(P)^{*^l}$ are given by the images of $u$ and
$v$ in $\kappa(P)$.

\paragraph*{Proposition 1.3}([S3], 2.5) If the index of $\alpha$
is $l$, then there are no hot points for $\alpha$.

\vskip 5mm

Suppose $P$ is a chilli point. Then $v(P) = u(P)^s$ for some $s$
with $1 \leq s \leq l-1$ and $s$ is called the {\it coefficient}
of $P$ with respect to $\pi$. To get some compatibility for these
coefficients, Saltman associates to $\alpha$ and $\calx$ the
following graph: The set of vertices is the set of irreducible
curves in ram$_{\calx}(\alpha)$ and there is an edge between two
vertices if there is a chilli point in the intersection of the two
irreducible curves corresponding to the vertices. A loop in this
graph is called a {\it chilli loop}.

\paragraph*{Proposition 1.4} ([S3], 2.6, 2.9) There exists a
projective model $\calx$ of $K$ such that there are no chilli
loops and no cool points on $\calx$ for $\alpha$.

\vskip 5mm

Let $F$ be a field of characteristic not equal to 2. The
$u$-invariant of $F$, denoted by $u(F)$, is defined as follows:
$$
u(F) = {\rm sup}\{\mbox{rk}(q) \mid q {\rm ~an  ~anisotropic~
quadratic ~form ~over}~ F \}.
$$
For $a_1, \cdots ,a_n \in F^*$, we denote the diagonal quadratic
form $a_1X_1^2 +  \cdots + a_nX_n^2$ by $<\! a_1, \cdots ,a_n
\!>$. Let $W(F)$ be the Witt ring of quadratic forms over $F$ and
$I(F)$ be the ideal of $W(F)$ consisting of even dimension forms.
Let $I^n(F)$ be the $n^{\rm th}$ power of the ideal $I(F)$. For
$a_1, \cdots, a_n \in F^*$, let $<\!<\! a_1, \cdots ,a_n\!>\!>$
denote the n-fold Pfister form $ <\! 1, a_1\!>\otimes \cdots
\otimes <\!1, a_n\!>$.  The abelian group $I^n(F)$ is generated by
$n$-fold Pfister forms.  The dimension modulo 2 gives an
isomorphism $e_0 : W(F)/I(F) \to H^0(F,
{\mathbf{Z}}/2{\mathbf{Z}})$. The discriminant gives an
isomorphism $e_1 : I(F)/I^2(F) \to H^1(F,
{\mathbf{Z}}/2{\mathbf{Z}})$. The classical result of Merkurjev
([M]), asserts that the Clifford invariant gives an isomorphism
$e_2 : I^2(F)/I^3(F) \to H^2(F, {\mathbf{Z}}/2{\mathbf{Z}})$.

Let $P_n(F)$ be the set of isometry classes of $n$-fold Pfister
forms over $F$. There is a well-defined map ([A])
$$
e_n : P_n(F) \to H^n(F, {\mathbf{Z}}/2{\mathbf{Z}})
$$
given by $e_n(<\! 1, a_1\! > \otimes \cdots <\! 1, a_n \!>) =
(-a_1) \cdots (-a_n) \in H^n(F,{\mathbf{Z}}/2{\mathbf{Z}} )$.

A quadratic form version of the Milnor conjecture asserts that
$e_n$ induces a surjective homomorphism $I^n(F) \to
H^n(F,{\mathbf{Z}}/2{\mathbf{Z}} )$ with kernel $I^{n+1}(F)$. This
conjecture was proved by Voevodsky, Orlov and Vishik. In this
paper we are interested in fields of 2-cohomological dimension  at
most 3. For such fields, the above Milnor's conjecture is already
proved by Arason, Elman and Jacob ([AEJ], Corollary 4 and Theorem
2), using the theorem of Markurjev ([M]).

Let $q_1$ and $q_2$ be two quadratic forms over $F$. We write $q_1
= q_2$ if they represent the same class in the Witt group $W(F)$.
We write $q_1 \simeq q_2$, if $q_1$ and $q_2$ are isometric
quadratic forms. We note that if the dimensions of $q_1$ and $q_2$
are equal and $q_1 = q_2$, then $q_1 \simeq q_2$.

\section*{2. Divisors on Arithmetic Surfaces}
In this section we recall a few results from a paper of Saltman
([S3]) on  divisors on arithmetic surfaces.

Let $\calz$ be a  connected, reduced scheme of finite type over a
Noetherian ring. Let $\calo^*_{\calz}$ be the sheaf of units in
the structure sheaf $\calo_{\calz}$. Let $\calp$ be a finite set
of closed points of $\calz$. For each $P \in \calp$, let
$\kappa(P)$ be the residue field at $P$ and $\iota_P : {\rm
Spec}(\kappa(P)) \to \calz$ be the natural morphism. Consider the
sheaf  $\calp^* = \oplus_{P \in \calp}\iota_P^*\kappa(P)^*$, where
$\kappa(P)^*$ denotes the group of units in $\kappa(P)$. Then
there is a surjective morphism of sheaves $\calo^*_{\calz} \to
\calp^*$ given by the evaluation at each $P \in \calp$. Let
$\calo_{\calz, \calp}^{*(1)}$ be its kernel. When there is no
ambiguity we denote $\calo_{\calz,\calp}^{*(1)}$ by
$\calo_{\calp}^{*(1)}$. Let $\calk$ be the sheaf of total quotient
rings on $\calz$ and $\calk^*$ be the sheaf of groups given by
units in $\calk$. Every element $\gamma \in H^0(\calz,
\calk^*/\calo^*)$ can be represented by a family $\{U_i, f_i \}$,
where $U_i$ are open sets in $\calz$, $f_i \in \calk^*(U_i)$ and
$f_if_j^{-1}  \in \calo^*(U_i \cap U_j)$. We say that an element
$\gamma = \{U_i, f_i \}$ of $H^0(\calz, \calk^*/\calo^*)$ {\it
avoids} $\calp$ if each $f_i$ is a unit at $P$ for all $P \in U_i
\cap \calp$. Let $H^0_{\calp}(\calz, \calk^*/\calo^*)$ be the
subgroup of $H^0(\calz, \calk^*/\calo^*)$ consisting of those
$\gamma$ which avoid $\calp$. Let $K^* = H^0(\calz, \calk^*)$ and
$K^*_{\calp}$ be the subgroup of $K^*$ consisting of those
functions which are units at all $P \in \calp$. We have a natural
inclusion $K^*_{\calp} \to H^0_{\calp}(\calz, \calk^*/\calo^*)
\oplus (\displaystyle{\oplus_{P \in \calp}}\kappa(P)^*)$.

We have the following

\paragraph*{Proposition 2.1} ([S3], 1.6) Let $\calz$ be a
connected, reduced scheme of finite type over a Noetherian ring.
Then $$H^1(\calz, \calo_{\calp}^{*(1)}) \simeq
\frac{H^0_{\calp}(\calz, \calk^*/\calo^*) \oplus (\oplus_{P \in
\calp}\kappa(P)^*)}{K^*_{\calp}}.$$

Let $k$ be a $p$-adic field and ${\cal O}_k$ the ring of integers
of $k$. Let $\calx$ be a connected regular surface with a
projective morphism $\eta: \calx \to {\rm Spec}({\cal O}_k)$. Let
$\overline{\calx}$ be the reduced special fibre of $\eta$. Assume
that $\overline{\calx}$ is connected. Note that $\overline{\calx}$
is connected if the generic fibre is geometrically integral. Let
$\calp$ be a finite set of closed points in $\calx$. Since every
closed point of $\calx$ is in $\overline{\calx}$, $\calp$ is also
a subset of closed points of $\overline{\calx}$. Let $m$ be an
integer coprime with $p$.

\paragraph*{Proposition 2.2} ([S3], 1.7) The canonical map
$H^1(\calx, \calo_{\calx, \calp}^{*(1)}) \to H^1(\overline{\calx},
\calo_{\overline{\calx}, \calp}^{*(1)})$ induces an isomorphism
$$ \frac{H^1(\calx, \calo_{\calx, \calp}^{*(1)})}{m H^1(\calx,
\calo_{\calx, \calp}^{*(1)})} \simeq \frac{H^1(\overline{\calx},
\calo_{\overline{\calx}, \calp}^{*(1)})}{m H^1(\overline{\calx},
\calo_{\overline{\calx}, \calp}^{*(1)})}.
$$

Let $\calx$ be as above. Suppose that $\overline{\calx}$ is
a union of regular curves $F_1, \cdots ,F_m$ on $\calx$ with only
normal crossings. Let $\calp$ be a finite set of closed points of
$\calx$ including all the points of $F_i \cap F_j$, $i \neq j$ and
at least one point from each $F_i$. Let $E$ be a divisor on
$\calx$ whose support does not pass through any point of $\calp$.
In particular, no $F_i$ is in the support of $E$. Hence there are
only finitely many closed points $Q_1, \cdots ,Q_n$ on the support
of $E$. For each closed point $Q_i$ on the support of $E$, let
$D_i$ be a regular geometric curve (i.e. a curve not contained in
the special fibre) on ${\calx}$ such that $Q_i$ is the
multiplicity one intersection of $D_i$ and $\overline{\calx}$.
Such a curve exists by ([S3], 1.1). We note that any closed point
on $\calx$ is a point of codimension 2 and there is a unique
closed point on any geometric curve on ${\calx}$ (cf. \S 1).

The following  is extracted from ([S3], \S 5).

\paragraph*{Proposition 2.3} Let $\calx, \calp, E, Q_i, D_i$ be as
above. For each closed point $Q_i$, let $m_i$ be the intersection
multiplicity of the support of $E$ and the special fibre
$\overline{\calx}$ at $Q_i$.  Let $l$ be a prime not equal to $p$.
Then there exist $\nu \in K^*$ and a divisor $E'$ on $\calx$ such
that
$$
(\nu) = -E + \sum_{i = 1}^n m_iD_i + lE'
$$
and $\nu(P) \in \kappa(P)^{*^l}$ for each $P \in \calp$.

\paragraph*{Proof.} Let $F$ be the divisor on $\calx$ given by
$\sum F_i$. Let $\gamma \in Pic(\calx)$ be the line bundle
equivalent to the class of the divisor $-E$ and $\overline{\gamma}
\in Pic(\overline{\calx})$ its image. Since the support of $E$
does not pass through the points of $\calp$ and $\calp$ contains
all the points of intersection of distinct $F_i$, $E$ and $F$
intersect only at smooth points of $\overline{\calx}$. In
particular, we have $\overline{\gamma} = -\sum m_i Q_i$. Let
$\gamma' \in H^1(\calx, \calo^*_{\calp})$ be the element which,
under the isomorphism of (2.1), corresponds to the class of the
element $(-E + \sum m_iD_i, 1)$ in $H^0_{\calp}(\calx,
\calk^*/\calo^*) \oplus (\oplus_{P \in \calp}\kappa(P)^*)$. Since
the $m_i$'s are intersection multiplicities of $E$ and
$\overline{\calx}$ at $Q_i$ and the image of $\sum m_iD_i$ in
$H^0_{\calp}(\overline{\calx}, \calk^*/\calo^*)$ is $\sum m_iQ_i$,
the image $\overline{\gamma'}$ of $\gamma'$ in
$H^1(\overline{\calx}, \calo^*_{\calp})$ is zero. By (2.2), we
have $\gamma' \in lH^1(\calx, \calo^*_{\calp})$. Using (2.1),
there exists $(E', (\lambda_P)) \in H^0_{\calp}(\calz,
\calk^*/\calo^*) \oplus (\oplus_{P \in \calp}\kappa(P)^*)$ such
that $(-E + \sum m_iD_i, 1) = l (E', (\lambda_P)) = (lE',
(\lambda_P^l))$ modulo $K^*_{\calp}$.  Thus there exists $\nu \in
K^*_{\calp} \subset K^*$ such that $(\nu) = (-E + \sum m_iD_i, 1)
- (lE', (\lambda_P^l))$. i.e. $(\nu) =  -E + \sum m_iD_i - lE'$
and $\nu(P) = \lambda_P^l$ for each $P \in \calp$. \hfill $\Box$

\section*{3. A local-global principle}
Let $k$ be a $p$-adic field, ${\cal O}_k$ be its ring of integers
and $K$ the function field of a  smooth, projective, geometrically
integral curve over $k$. Let $l$ be a prime not equal to $p$.
Throughout this section, except in 3.6, we assume that $k$
contains a primitive $l^{\rm th}$ root of unity. We fix a
generator $\rho$ for $\mu_l$ and identify $\mu_l$ with
${\mathbf{Z}}/l{\mathbf{Z}}$.

\paragraph*{Lemma  3.1} Let $\alpha \in H^2(K, \mu_l)$. Let
${\cal X}$ be a regular projective model of $K$. Assume that the
ramification locus $\rama$ is a union of regular curves $\{ C_1,
\cdots ,C_r\}$ with only normal crossings. Let $T$ be a finite set
of closed points of ${\cal X}$ including the points of $C_i \cap
C_j$, for all $i \neq j$. Let $D$ be an irreducible curve on
${\cal X}$ which is not in the ramification locus of $\alpha$ and
does not pass through any point in $T$. Then $D$ intersects $C_i$
at  points $P$  where $\pa$ is unramified. Suppose further that at
such points $P$, $\partial_{C_i}(\alpha)$ specializes to $0$ in
$H^1(\kappa(P), {\mathbf{Z}}/l{\mathbf{Z}})$. Then $\alpha$ is
unramified at $D$ and specializes to 0 in $H^2(\kappa(D), \mu_l)$.

\paragraph*{Proof.} Since $k$ contains a primitive $l^{\rm th}$
of unity, we fix a generator $\rho$ for $\mu_l$ and identify the
Galois modules $\mu_l^{\otimes j}$ with
${\mathbf{Z}}/l{\mathbf{Z}}$. Let $P$ be a point in the
intersection of $D$ and the support of ${\rm
ram}_{\calx}(\alpha)$. Since $D$ does not pass through the points
of $T$ and $T$ contains all the points of intersection of distinct
$C_j$, the point $P$ belongs to a unique curve $C_i$ in the
support of ram$_{\calx}(\alpha)$. We have ([S1], 1.2) $\alpha =
\alpha' + (u, \pi)$, where $\alpha'$ is unramified on ${\cal
O}_{{\cal X}, P}$, $u \in {\cal O}_{{\cal X}, P}$ is a unit and
$\pi \in {\cal O}_{{\cal X}, P}$ is a prime defining the curve
$C_i$ at $P$. Therefore $\pa = (\overline{u}) \in H^1(\kappa(C_i),
{\mathbf{Z}}/l{\mathbf{Z}})$ is unramified at $P$.

Suppose that $\partial_{C_i}(\alpha)$ specialises to zero in
$H^1(\kappa(P), {\mathbf{Z}}/l{\mathbf{Z}})$. Since $D$ is not in
the ramification locus of $\alpha$, $\alpha$ is unramified at $D$.
Let $\overline{\alpha}$ be the specialization of $\alpha$ in
$H^2(\kappa(D), {\mathbf{Z}}/l{\mathbf{Z}})$. Since $\kappa(D)$ is
either a $p$-adic field or a function field of a curve over a
finite field, to show that $\overline{\alpha}$ is zero, by class
field theory it is enough to show that $\overline{\alpha}$ is
unramified at every discrete valuation of $\kappa(D)$.

Let $v$ be a discrete valuation of $\kappa(D)$ and $R$ the
corresponding discrete valuation ring. Then there exists a closed
point $P$ of $D$ such that $R$ is a localization of the integral
closure of the one dimensional local ring ${\cal O}_{D,P}$ of $P$
on $D$. The local ring ${\cal O}_{D,P}$ is a quotient of the local
ring ${\cal O}_{{\cal X}, P}$.

Suppose $P$ is not on the ramification locus of $\alpha$. Then
$\alpha$ is unramified on ${\cal O}_{{\cal X}, P}$ and hence
$\overline{\alpha}$ on  $\overline{\cal O}_{D,P}$. In particular
$\overline{\alpha}$ is unramified at $R$.

Suppose $P$ is on the ramification locus of $\alpha$. As before,
we have  $\alpha = \alpha' + (u, \pi)$, where $\alpha'$ is
unramified on ${\cal O}_{{\cal X}, P}$, $u \in {\cal O}_{{\cal X},
P}$ is a unit and $\pi \in {\cal O}_{{\cal X}, P}$ is a prime
defining the curve $C_i$ at $P$. Therefore $\pa = \overline{u}$ in
$\kappa(C_i)^*/\kappa(C_i)^{*^l}$. Since, by the assumption, $\pa$
specializes to $0$  at $P$, $u(P) \in \kappa(P)^{*^l}$. We have
$\overline{\alpha} = \overline{\alpha'} + (\overline{u},
\overline{\pi}) \in H^2(\kappa(D), {\mathbf{Z}}/l{\mathbf{Z}})$.
Since $\alpha'$ is unramified at $P$, the residue of
$\overline{\alpha}$ at $R$ is $(u(P))^{\nu(\overline{\pi})}$.
Since $\kappa(P)$ is contained  in the residue field of the
discrete valuation ring $R$ and $u(P)$ is an $l^{th}$ power in
$\kappa(P)$, it follows that $\overline{\alpha}$ is unramified at
$R$. \hfill $\Box$

\paragraph*{Proposition 3.2} Let $K$ and $l$ be as above. Let $\alpha
\in H^2(K, \mu_l)$ with index $l$. Let ${\cal X}$ be a regular
projective model of $K$ such that the ramification locus $\rama$
and the special fibre of ${\cal X}$ are a union of regular curves
with only normal crossings  and $\alpha$ has no cool points and no
chilli loops on $\calx$ (cf. 1.4).  Let $s_i$ be the corresponding
coefficients (cf. \S 1). Let $F_1, \cdots , F_r$ be irreducible
regular curves on ${\cal X}$ which are not in $\rama = \{ C_1,
\cdots ,C_n\}$ and such that $\{ F_1, \cdots ,F_r \} \cup \rama$
have only normal crossings. Let $m_1, \cdots, m_r$ be integers.
Then there exists $f \in K^*$ such that
$$
{\rm div}_{\cal X}(f) = \sum s_iC_i + \sum m_{s}F_{s} + \sum
n_jD_j + lE',
$$
where $D_1, \cdots , D_t$ are irreducible curves which are not
equal to $C_i$ and $F_s$ for all $i$ and $s$ and $\alpha$
specializes to zero at $D_j$ for all $j$  and $(n_j, l) = 1$.

\paragraph*{Proof.} Let $T$ be a finite set of closed points of $\calx$
containing all the points of intersection of distinct $C_i$ and
$F_s$ and at least one point from each $C_i$ and $F_s$. By a
semilocal argument, we choose $g \in K^*$ such that ${\rm
div}_{\cal X}(g) = \sum s_iC_i + \sum m_sF_s + G$ where $G$ is a
divisor on ${\cal X}$ whose support does not contain any of $C_i$
or $F_s$ and does not intersect $T$.

Since $\alpha$ has no cool points and no chilli loops on $\calx$,
by ([S3], Prop. 4.5), there exists $u \in K^*$ such that ${\rm
div}_{\cal X}(ug) = \sum s_iC_i + \sum m_sF_s + E$, where $E$ is a
divisor on ${\cal X}$ whose support does not contain any $C_i$ or
$F_s$, does not pass through the points in $T$ and either $E$
intersect $C_i$ at a point $P$ where the specialization of $\pa$
is $0$ or the intersection multiplicity $(E \cdot C_i)_P$ is a
multiple of $l$.

Suppose $C_i$ for some $i$ is a geometric curve on $\calx$. Since
every closed point of $\calx$ is on the special fibre
$\overline{\calx}$, the closed point of $C_i$ is in $T$. Since the
support of $E$ avoids all the points in $T$, the support of $E$
does not intersect $C_i$. Thus the support of $E$ intersects only
those $C_i$ which are in the special fibre $\overline{\calx}$. Let
$Q_1, \cdots , Q_t$ be the points of intersection of the support
of the divisor $E$ and the special fibre with intersection
multiplicity $n_j$ at $Q_j$ coprime with $l$, for $1 \leq j \leq
r$. In particular we have $\partial_{C_j}(\alpha)  = 0$ for $1
\leq j \leq r$. For each $Q_j$, let $D_j$ be a regular geometric
curve on $\calx$ such that $Q_j$ is the multiplicity one
intersection of $D_j$ and $\overline{\calx}$ (cf. paragraph after
2.2). Then by (2.3) there exists $\nu \in K^*$ such that ${\rm
div}_{\cal X}(\nu) = -E + \sum n_jD_j + lE'$ and $\nu(P) \in
\kappa(P)^{*^l}$ for all $P \in T$. Let $f = ug\nu \in K^*$. Then
we have
$$
{\rm div}_{\cal X}(f) = \sum s_iC_i + \sum m_{s}F_{s} + \sum
n_jD_j + lE'.
$$

Since each $Q_j$ is the only closed point on $D_j$ and $\pa$
specialises to zero at $Q_j$, by (3.1), the $\alpha$ specializes
to 0 at $D_j$. Thus $f$ has all the required properties. \hfill
$\Box$

\paragraph*{Lemma 3.3} Let $\alpha \in H^2(K,  \mu_l)$
and let $v$ be a discrete valuation of $K$. Let $u \in K^*$ be a
unit at $v$ such that  $\overline{u} \in \kappa(v)^* \setminus
\kappa(v)^{*^l}$. Suppose further that if  $\alpha$ is
ramified at $v$,  $\partial_v(\alpha) = [L] \in H^1(\kappa(v),
{\mathbf{Z}}/l{\mathbf{Z}})$, where $L = K(u^{\frac{1}{l}})$.
Then, for any  $g \in L^*$, the image of $(N_{L/K}(g)) \cdot
\alpha \in H^3(K_v, \mu_l^{\otimes 2})$ is zero.

\paragraph*{Proof.} We identify the Galois modules $\mu_l^{\otimes j}$
with ${\mathbf{Z}}/l{\mathbf{Z}}$ as before. Since $u$ is a unit
at $v$ and $\overline{u} \not\in \kappa(v)^{*^l}$, there is a
unique discrete valuation $\tilde{v}$ of $L$ extending the
valuation $v$ of $K$, which is unramified with residual degree
$l$. In particular $v(N_{L/K}(g))$ is a multiple of $l$. Thus if
$\alpha' \in H^2(K_v, {\mathbf{Z}}/l{\mathbf{Z}})$ is unramified
at $v$, then $(N_{L/K}(g)) \cdot \alpha' \in H^3(K_v,
{\mathbf{Z}}/l{\mathbf{Z}})$ is unramified. Since $H^3_{nr}(K_v,
{\mathbf{Z}}/l{\mathbf{Z}}) = 0$ (cf. 1.2), we have $(N_{L/K}(g))
\cdot \alpha' = 0$ for any $\alpha' \in H^2(K_v,
{\mathbf{Z}}/l{\mathbf{Z}})$ which is unramified at $v$. In
particular, if $\alpha$ is unramified at $v$, then $\alpha \cdot
(N_{L/K}(g)) = 0$.

Suppose that $\alpha$ is ramified at $v$. Then by the choice of
$u$, we have $\alpha = \alpha' + (u) \cdot (\pi_v)$, where $\pi_v$
is a parameter at $v$ and $\alpha' \in H^2(K_v,
{\mathbf{Z}}/l{\mathbf{Z}})$ is unramified at $v$. Thus we have
$(N_{L/K}(g)) \cdot \alpha = (N_{L/K}(g)) \cdot \alpha' +
(N_{L/K}(g)) \cdot (u) \cdot (\pi_v) = (N_{L/K}(g)) \cdot (u)
\cdot (\pi_v) \in H^3(K_v, {\mathbf{Z}}/l{\mathbf{Z}})$. Since
$L_v = K_v(u^{\frac{1}{l}})$, we have $((N_{L/K}(g)) \cdot (u)  =
0 \in H^2(K_v, {\mathbf{Z}}/l{\mathbf{Z}})$ and $(N_{L/K}(g))
\cdot \alpha = 0$ in $H^3(K_v, {\mathbf{Z}}/l{\mathbf{Z}})$.
\hfill $\Box$

\paragraph*{Theorem 3.4} Let $K$ and $l$ be as above. Let $\alpha
\in H^2(K, \mu_l)$ and $\zeta \in H^3(K,
 \mu_l^{\otimes 2})$. Assume that the index of $\alpha$ is
$l$. Let $\calx$ be a regular projective model of $K$. Suppose
that for each $x \in \calx^1$, there exists $f_x \in K_x^*$ such
that $\zeta  = \alpha  \cup (f_x) \in H^3(K_x,
 \mu_l^{\otimes 2})$, where $K_x$ is the completion of $K$
at the discrete valuation given by  $x$. Then there exists $f \in
K^*$ such that $\zeta = \alpha \cup (f) \in H^3(K, \mu_l^{\otimes
2})$.

\paragraph*{Proof.}  We identify the Galois modules $\mu_l^{\otimes j}$
with ${\mathbf{Z}}/l{\mathbf{Z}}$ as before. By  weak
approximation, we may find $f \in K^*$ such that $(f) = (f_v) \in
H^1(K_v, {\mathbf{Z}}/l{\mathbf{Z}})$ for all the discrete
valuations corresponding to the irreducible curves in $ram_{\cal
X}(\alpha) \cup ram_{\cal X}(\zeta)$. Let
$${\rm
div}_{\cal X}(f) = C' + \sum m_iF_i + lE,
$$
where $C'$ is a  divisor with support contained in $ram_{\cal
X}(\alpha) \cup ram_{\cal X}(\zeta)$, $F_i$'s are distinct
irreducible curves which are not in $ram_{\cal X}(\alpha) \cup
ram_{\cal X}(\zeta)$, $m_i$ coprime with $l$ and $E$ some divisor
on ${\cal X}$.

For any $C_j \in ram_{\cal X}(\zeta) \setminus ram_{\cal
X}(\alpha)$, let $\lambda_j \in \kappa(C_j)^* \setminus
\kappa(C_j)^{*^l}$. By weak approximation, we choose $u \in K^*$
with $\overline{u}= \pa \in
H^1(\kappa(C_i),{\mathbf{Z}}/l{\mathbf{Z}}) $ for all $C_i \in
ram_{\cal X}(\alpha)$, $\nu_{F_i}(u) = m_i$, where $\nu_{F_i}$ is
the discrete valuation at $F_i$ and $\overline{u} = \lambda_j$ for
any $C_j \in ram_{\cal X}(\zeta) \setminus ram_{\cal X}(\alpha)$.
In particular $u$ is a unit at the generic point of $C_j$ and
$\overline{u} \not\in \kappa(C_j)^{*^l}$ for any $C_j \in
ram_{\cal X}(\zeta) \setminus ram_{\cal X}(\alpha)$.

Let $L = K(u^{\frac{1}{l}})$. Let $\eta : {\cal Y} \to \calx$ be
the normalization  of ${\cal X}$ in $L$. Since $\nu_{F_i} (u) =
m_i$ and $m_i$ is coprime with $l$, $\eta: {\cal Y} \to \calx$ is
ramified at $F_i$. In particular  there is a unique irreducible
curve $\tilde{F_i}$ in ${\cal Y}$ such that $\eta(\tilde{F_i}) =
F_i$ and  $\kappa(F_i) = \kappa(\tilde{F_i})$.

Let $\pi : \tilde{\cal Y} \to {\cal Y}$ be a proper birational
morphism such that the ramification locus $ram_{\tilde{\cal
Y}}(\alpha_L)$ of $\alpha_L$ on $\tilde{\cal Y}$ and the strict
transform of the curves $\tilde{F_i}$ on $\tilde{\cal Y}$ is a
union of regular curves with only normal crossings and there are
no cool points and no chilli loops for $\alpha_L$ on $\tilde{\cal
Y}$ (cf 1.4). We denote the strict transforms of $\tilde{F_i}$ by
$\tilde{F_i}$ again. By (3.2), there exists $g \in L^*$ such that
$$
div_{\tilde{\cal Y}}(g) = C + \sum -m_i\tilde{F_i} + \sum n_jD_j +
lD,
$$
where the support of $C$ is contained in $ram_{\tilde{\cal
Y}}(\alpha_L)$ and $D_j$'s are irreducible curves which are not in
$ram_{\tilde{\cal Y}}(\alpha_L)$ and $\alpha_L$ specializes to
zero at all $D_j$'s.

We now claim that $\zeta = \alpha \cup (fN_{L/K}(g))$. Since the
group $H^3_{nr}(K/\calx, {\mathbf{Z}}/l{\mathbf{Z}}) = 0$ ([K],
5.2), it is enough to show that $\zeta - \alpha \cup
(fN_{L/K}(g))$ is unramified on ${\cal X}$. Let $S$ be an
irreducible curve on ${\cal X}$. Since the residue map
$\partial_S$ factors through the completion $K_S$, it suffices to
show that $\zeta - \alpha \cup (fN_{L/K}(g)) = 0$ over $K_S$.

Suppose $S$ is not in $ram_{\cal X}(\alpha) \cup ram_{\cal
X}(\zeta) \cup Supp(fN_{L/K}(g))$. Then each of $\zeta$ and
$\alpha \cdot(fN_{L/K}(g))$ is unramified at $S$.

Suppose that $S$ is in $ram_{\cal X}(\alpha) \cup ram_{\cal
X}(\zeta)$. Then by the choice of $f$ we have $(f) = (f_v) \in
H^1(K_v, {\mathbf{Z}}/l{\mathbf{Z}})$ where $v$ is the discrete
valuation associated to $S$. Hence $\zeta = \alpha \cup (f)$ over
the completion  $K_S$ of $K$ at the discrete valuation given by
$S$. It follows from (3.3) that $(N_{L/K}(g)) \cup \alpha = 0$
over $K_S$ and $\zeta = \alpha \cup (f N_{L/K}(g))$ over $K_S$.

Suppose that $S$ is in the support of $div_{\cal X}(fN_{L/K}(g))$
and not in $ram_{\cal X}(\alpha) \cup ram_{\cal X}(\zeta)$. Then
$\alpha$ is unramified at $S$. We show that in this case $\alpha
\cup (fN_{L/K}(g)) = 0$ over $K_S$. We have
$$
\begin{array}{lll}
div_{\cal X}(fN_{L/K}(g)) & = &  div_{\cal X}(f) +   div_{\cal
X}(N_{L/K}(g)) \cr & = & C' + \sum m_iF_i + lE + \eta_* \pi_*(C +
\sum -m_i\tilde{F_i} + \sum n_jD_j + lD)  \cr & = & C' +
\eta_*\pi_*(C)  + \sum n_j \eta_*\pi_*(D_j) + l \eta_*\pi_*(E).
\end{array}
$$
We note that if $D_j$ maps to a point, then $\eta_*\pi_*(D_j) =
0$. Since the support of $C$ is contained in $ram_{\tilde{\cal
Y}}(\alpha_L)$, the support of $\eta_*\pi_*(C)$ is contained in
$ram_{\cal X}(\alpha)$. Thus $S$ is in the support of
$\eta_*\pi_*(D_j)$ for some $j$ or $S$ is in the support of $l
\eta_*\pi_*(E)$. In the later case, clearly $\alpha \cdot
(fN_{L/K}(g))$ is unramified at $S$ and hence $\alpha \cdot
(fN_{L/K}(g)) = 0$ over $K_S$. Suppose $S$ is in the support of $
\eta_*\pi_*(D_j)$ for some $j$. In this case, if $D_j$ lies over
an inert curve, then $\eta_*\pi_*(D_j)$ is a multiple of $l$ and
we are done. Suppose that $D_j$ lies over a split curve. Since
$\alpha_L$ specializes to zero at $D_j$, it follows that $\alpha$
specializes to zero at $\eta_*\pi_*(D_j)$ and we are done. \hfill
$\Box$

\paragraph*{Theorem 3.5.} Let $k$ be a $p$-adic field and $K$ a function
field of a curve over $k$. Let $l$ be a prime not equal to $p$.
Suppose that all the $l^{\rm th}$ roots of unity are in $K$. Then
every element in $H^3(K, \mu_l^{\otimes 3})$ is a symbol.

\paragraph*{Proof.}   We again identify the Galois modules
$\mu_l^{\otimes j}$ with ${\mathbf{Z}}/l{\mathbf{Z}}$.

Let $v$ be a discrete valuation of $K$ and $K_v$ the completion of
$K$ at $v$. By (1.2), every element in $H^3(K_v,
{\mathbf{Z}}/l{\mathbf{Z}})$ is a symbol.

Let $\zeta \in H^3(K, {\mathbf{Z}}/l{\mathbf{Z}})$ and  $\calx$ be
a regular projective model of $K$.  Let $v$ be a discrete
valuation of $K$ corresponding to an irreducible curve in
ram$_{\calx}(\zeta)$. Then we have $\zeta = (f_v) \cdot (g_v)
\cdot (h_v)$ for some $f_v, g_v, h_v \in K^*_v$. By weak
approximation, we can find $f, g \in K^*$ such that $(f) = (f_v)$
and $(g) = (g_v)$ in $H^1(K_v, {\mathbf{Z}}/l{\mathbf{Z}})$ for
all discrete valuations $v$ corresponding to the irreducible
curves in ram$_{\calx}(\zeta)$. Let $v$ be a discrete valuation of
$K$ corresponding to an irreducible curve $C$ in $\calx$. If $C$
is in ram$_{\calx}(\zeta)$, then by the choice of $f$ and $g$ we
have $\zeta = (f) \cdot (g) \cdot (h_v) \in H^3(K_v,
{\mathbf{Z}}/l{\mathbf{Z}})$. If $C$ is not in the
ram$_{\calx}(\zeta)$, then $\zeta \in H^3_{\rm
nr}(K_v,{\mathbf{Z}}/l{\mathbf{Z}}) \simeq H^3(\kappa(v),
{\mathbf{Z}}/l{\mathbf{Z}}) = 0$. In particular we have $\zeta =
(f) \cdot (g) \cdot (1) \in H^3(K_v, {\mathbf{Z}}/l{\mathbf{Z}})$.
Let $\alpha = (f) \cdot (g) \in H^2(K,
{\mathbf{Z}}/l{\mathbf{Z}})$. Then we have $\zeta = \alpha \cdot
(h'_v) \in H^3(K_v, {\mathbf{Z}}/l{\mathbf{Z}})$ for some $h'_v
\in K^*_v$ for each discrete valuation $v$ of $K$ associated to
any point of $\calx^1$. By (3.4), there exists $h \in K^*$ such
that $\zeta = \alpha \cdot (h) =  (f) \cdot (g) \cdot (h) \in
H^3(K, {\mathbf{Z}}/l{\mathbf{Z}}).$ \hfill $\Box$

\paragraph*{Remark 3.6.} We remark that all the results of this section
can be extended to the situation where $k$ does not necessarily
contain a primitive $l^{\rm th}$ root of unity. This can be
achieved by going to the extension $k'$ of $k$ obtained by
adjoining a primitive $l^{\rm th}$ of unity to $k$ and noting that
the extension $k'/k$ is unramified of degree $l-1$. We do not use
this remark in the sequel.

\section*{4. The $u$-invariant }
In (4.1) and (4.2) below, we give some necessary conditions for a
field $k$ to have the u-invariant less than or equal to 8. If $K$
is the function field of a curve over a $p$-adic field and $K_v$
is the completion of $K$ at a discrete valuation $v$ of $K$, then
the residue field $\kappa(v)$ of $K_v$, which is either a global
field of positive characteristic or a $p$-adic field, has
$u$-invariant 4. By a theorem of Springer, $u(K_v) = 8$ and we use
(4.1) and (4.2)  for $K_v$.

\paragraph*{Proposition 4.1} Let $K$ be a field of characteristic not equal to
2. Suppose that $u(K) \leq 8$. Then $I^4(K) = 0$ and  every
element in $I^3(K)$ is a 3-fold Pfister form. Further if $\phi$ is
a 3-fold Pfister form  and $q_2$  a rank 2 quadratic form over
$K$, then there exists $f, g, h \in K^*$ such that $f$ is a value
of $q_2$ and $\phi  = <\! 1, f\!><\! 1, g\!><\! 1,h\!>  $.

\paragraph*{Proof.} Suppose that $u(K) = 8$. Then every 4-fold Pfister form is
isotropic and hence hyperbolic; in particular, $I^4(K) = 0$. Let
$\phi_1 = <\! 1, f_1\!><\! 1, g_1\!><\! 1,h_1\!>$ and $\phi_2 =
<\!1, f_2\!><\! 1, g_2\!><1,h_2\!>$ be two anisotropic 3-fold
Pfister forms. Since $u(K) \leq 8$, the Witt index of $\phi_1 -
\phi_2$ is at least 4. In particular $\phi_1 - \phi_2$ is
isotropic, i.e. $\phi_1$ and $\phi_2$ represent a common value $a
\in K$. Since $\phi_1$ and $\phi_2$ are anisotropic, $a \neq 0$.
We have $\phi_1 = <\! a \!> \perp \phi_1'$ and $\phi_2 = <\! a \!>
\perp \phi_2'$ for some quadratic form $\phi_1$ and $\phi_2$ over
$K$ ([Sc], p.7, Lemma 3.4). Since the Witt index of $\phi_1 -
\phi_2$ is at least 4, the Witt index of $\phi_1' - \phi_2'$ is at
least 3. Repeating this process, we see that  there exists a
quadratic form $<\! a,~b, ~c, ~d \!>$ over $K$ which is a subform
of both $\phi_1$ and $\phi_2$ . Let $C$ be a conic given by the
quadratic form $abc<\! a, b, c\!> = <\! bc, ac, ab\!>$. Since
$abc<\! a, b, c\!>$ is isotropic over the function field $K(C)$ of
$C$ (cf. [Sc], p.154, Remark 5.2(iv)) and $abc<\! a, b, c\!>$ is a
subform of $abc\phi_1$ and $abc\phi_2$, $\phi_1$ and $\phi_2$ are
isotropic over $K(C)$. Since $\phi_1$ and $\phi_2$ are Pfister
forms, they are hyperbolic over $k(C)$ ([Sc], p.144, Cor. 1.5).
Let $\psi = <\! 1, bc \!><\! 1, ac\!>$ and $k(\psi)$ be the
function field of the quadratic form $\psi$. Then $k(C)$ and
$k(\psi)$ are birational (cf. [Sc], p.154) and hence $\phi_1$ and
$\phi_2$ are hyperbolic over $k(\psi)$. Therefore $\phi_1 \simeq
<\! 1, bc\!><\! 1, ac\!><\! a_1, b_1\!>$ and $\phi_2 \simeq <\! 1,
bc\!><\! 1,ac\!><\! a_2, b_2\!>$ for some $a_1, a_2, b_1, b_2 \in
K^*$ (cf. [Sc], p.155, Th.5.4). Since $I^4(K) = 0$, we have
$\phi_1 \simeq <\! 1, bc\!><\! 1, ac\!><\! 1, a_1b_1\!>$ and
$\phi_2 \simeq <\! 1, bc\!><\!1, ac\!><\! 1, a_2b_2\!>$. We have
$$
\begin{array}{lll}
\phi_1 + \phi_2 & = & \phi_1 - \phi_2 \cr & = & <\! 1, bc\!><\! 1,
ac\!><\! 1, a_1b_1, -1, -a_2b_2\!> \cr & = &  <\! 1, bc\!><\! 1,
ac\!><\! a_1b_1, -a_2b_2\!> \cr & = & a_1b_1<\! 1, bc\!><\!
1,ac><\! 1, -a_1a_2b_1b_2\!> \cr &= &<\! 1, bc\!><\! 1, ac\!><\!
1, -a_1a_2b_1b_2\!>.
\end{array}$$
Thus the sum of any two 3-fold Pfister forms is a Pfister form in
$I^3(K)$ and every element in $I^3(K)$ is the class of  a 3-fold
Pfister form.

Let $\phi = <\! 1, a\!><\! 1, b\!><\! 1, c\!>$ be a 3-fold Pfister
form and $\phi'$ be its pure subform. Let $q_2$ be a quadratic
form over $K$ of dimension 2. Since  dim$(\phi') = 7$ and $u(K)
\leq 8$, the quadratic form $\phi' - q_2$ is isotropic. Therefore
there exists $f \in K^*$ which is a value of $q_2$ and $\phi'
\simeq <\! f\!> + \phi''$ for some quadratic form $\phi''$ over
$K$. Hence by ([Sc], p.143), $\phi = <\! 1, f\!><\! 1, b'\!><\! 1,
c'\!>$ for some $b', c' \in K^*$. \hfill $\Box$

\paragraph*{Proposition 4.2.} Let $K$ be a field of
characteristic not equal to 2. Suppose that $u(K) \leq 8$. Let
$\phi = <\! 1, f\!><\! 1, a\!><\! 1, b\!>$ be a 3-fold Pfister
form over $K$ and $q_3$ a quadratic form over $K$ of dimension 3.
Then there exist $g, h \in K^*$ such that $g$ is a value of $q_3$
and $\phi = <\! 1, f\!><\! 1, g\!><\! 1, h\!>$.

\paragraph*{Proof.} Let $\psi = <\! 1, f\!><\! a, b, ab\!>$.
Since $u(K) \leq 8$, the quadratic form $\psi - q_3$ is isotropic.
Hence there exists $g \in K^*$ which is a common value of $q_3$
and $\psi$. Thus, $\psi \simeq <\! g\!> + \psi_1$ for some
quadratic form $\psi_1$ over $K$. Since $\psi$ is hyperbolic over
$K(\sqrt{-f})$,  $\psi_1 \simeq <\! 1, f\!><\! a_1, b_1\!> + <\!
g_1\!>$ for some $a_1, b_1, g_1 \in K^*$. By comparing the
determinants, we get  $g_1 = gf$ modulo squares. Hence $\psi = <\!
1, f\!><\! g, a_1, b_1\!>$ and $\phi = <\! 1, f\!> + \psi = <\! 1,
f\!><\! 1, g, a_1, b_1\!>$. The form $\phi$ is isotropic and hence
hyperbolic over  the function field of the conic given by $<\! f,
g, fg\!>$. Hence, as in 4.1,   $\phi = \lambda<\! 1, f\!><\! 1,
g\!><\! 1, h\!>$ for some $\lambda, h \in K^*$. Since $I^4(K) =
0$, $\phi = <\! 1, f\!><\! 1, g\!><\! 1, h\!>$ with  $g$ a value
of $q_3$. \hfill $\Box$

\paragraph*{Proposition 4.3.} Let $K$ be a field of characteristic
not equal to 2. Assume the following:
\begin{enumerate}

\item{} Every element in $H^2(K, \mu_2)$ is a sum of at most 2
symbols.

\item{} Every element in $I^3(K)$ is  equal to a 3-fold Pfister
form.

\item{} If $\phi$ is a 3-fold Pfister form and $q_2$ is a
quadratic form over $K$ of dimension 2, then $\phi = <\! 1,
f\!><\! 1, g\!><\! 1, h\!>$ for some $f, g, h \in K^*$ with $f$ a
value of $q_2$.

\item{} If $\phi = <\! 1, f\!><\! 1,a\!><\! 1, b\!>$ is a 3-fold
Pfister form and $q_3$ a quadratic form over $K$ of dimension 3,
then $\phi = <\! 1, f\!><\! 1, g\!><\! 1, h\!>$ for some $g, h \in
K^*$ with $g$ a value of $q_3$.

\item{} $I^4(K) = 0.$
\end{enumerate}

Then $u(K) \leq 8$.

\paragraph*{Proof.} Let $q$ be a quadratic form over $K$ of dimension 9.
Since every element in $H^2(K, \mu_2)$ is a sum of at most 2
symbols, as in ([PS], proof of 4.5 ), we find a quadratic form
$q_5 = \lambda<\! 1, a_1, a_2, a_3, a_4\!>$ over $K$ such that
$\phi = q  +  q_5 \in I^3(K)$.  By the assumptions 2, 3 and 4,
there exist $f, g, h \in K^*$ such that $\phi = <\! 1, f\!><\! 1,
g\!><\! 1, h\!>$ and $f$ is a value of $<\! a_1, a_2\!>$ and $g$
is a value of $<\! fa_1a_2, a_3, a_4\!>$. We have $<\! a_1, a_2\!>
\simeq <\! f, fa_1a_2\!>$ and $<\! fa_1a_2, a_2, a_3\!> \simeq <\!
g, g_1, g_2\!>$ for some $g_1, g_2 \in K^*$. Since $I^4(K) = 0$,
we have $\lambda \phi = \phi$ and
$$
\begin{array}{rlll}
\lambda q & = & \lambda q + \lambda q_5 - \lambda q_5 \cr & = &
\lambda \phi - \lambda q_5 \cr & = & \phi - \lambda q_5 \cr & = &
<\! 1, f\!><\! 1, g\!><\! 1, h\!> - <\! 1, a_1, a_2, a_3, a_4\!>
\cr & = & <\! 1, f\!><\! 1, g\!><\! 1, h\!> - <\! 1, f, g, g_1,
g_2\!> \cr & = & <\! gf\!> + <\! 1,f\!><\! h, gh\!> - <\! g_1,
g_2\!>.
\end{array}
$$
Since the dimension of $\lambda q$ is 9 and the dimension of $<\!
gf\!> \perp <\! 1, f\!><\! h, gh\!> - <\! g_1, g_2\!>$ is 7, it
follows that $\lambda q$ and hence $q$ is isotropic over $K$.
\hfill $\Box$

\paragraph*{Proposition 4.4.} Let  $k$ be a $p$-adic field, $p \neq
2$ and $K$ a function field of a curve over $k$. Let $\phi$ be a
3-fold Pfister form over $K$ and $q_2$ a quadratic form over $K$
of dimension 2. Then there exist $f, a,  b \in K^*$ such that $f$
is a value of $q_2$ and $\phi  = <\! 1, f\!><\! 1, a\!><\! 1,
b\!>$.

\paragraph*{Proof.} Let $\zeta = e_3(\phi)  \in
H^3(K, \mu_2)$.  Let $\calx$ be a projective regular model of $K$.
Let $C$ be an irreducible curve on $\calx$ and $v$ the discrete
valuation given by $C$. Let $K_v$ be the completion of $K$ at $v$.
Since the residue field $\kappa(v) = \kappa(C)$ is either a
$p$-adic field or a function field of a curve over a finite field,
$u(\kappa(v)) = 4$ and  $u(K_v) = 8$ ([Sc], p.209). By (4.1),
there exist $f_v, a_v, b_v \in K_v^*$ such that $f_v$ is a value
of $q_2$ over $K_v$ and $\phi = <\! 1, f_v\!><\! 1, a_v\!><\! 1,
b_v \!>$ over $K_v$. By weak approximation, we can find $f,a \in
K^*$ such that $f$ is a value of $q_2$ over $K$ and $f = f_v, a =
a_v$ modulo $K_v^{*^2}$ for all discrete valuations $v$
corresponding to the irreducible curves $C$ in the support of
${\rm ram}_{\calx}(\zeta)$. Let $C$ be any irreducible curve on
${\calx}$ and $v$ the discrete valuation of $K$ given by $C$. If
$C$ is in the support of ${\rm ram}_{\calx}(\zeta)$, then by the
choice of $f$ and $a$, we have $\zeta = e_3(\phi) = (-f)\cdot (-a)
\cdot (-b_v)$ over $K_v$. If $C$ is not in the support of ${\rm
ram}_{\calx}(\zeta)$, then $\zeta \in H^3_{\rm nr}(K_v, \mu_2)
\simeq H^3(\kappa(v), \mu_2) = 0$. In particular we have $\zeta =
(-f) \cdot (-a) \cdot (1)$ over $K_v$. Let $\alpha = (-f) \cdot
(-a) \in H^2(K, \mu_2)$. By (3.4), there exists $b \in K^*$ such
that $\zeta = \alpha \cdot (-b) \in H^3(K, \mu_2)$. Since  $e_3 :
I^3(K) \to H^3(K, \mu_2)$ is an isomorphism, we have $ \phi = <\!
1, f\!><\! 1, a\!><\! 1, b\!>$ as required. \hfill $\Box$.

\vskip 5mm

There is a different proof of the Proposition 4.4 in ([PS], 4.4)!

\paragraph*{Proposition 4.5} Let  $k$ be a $p$-adic field, $p \neq
2$ and $K$ a function field of a curve over $k$. Let $\phi =  <\!
1, f\!><\! 1, a\!><\! 1, b\!>$ be a 3-fold Pfister form over $K$
and $q_3$ a quadratic form over $K$ of dimension 3. Then there
exist $g, h \in K^*$ such that $g$ is a value of $q_3$ and $\phi =
<\! 1, f\!><\! 1, g\!><\! 1, h\!>$.

\paragraph*{Proof.} Let $\zeta = e_3(\phi) = (-f) \cdot (-a) \cdot (-b)
\in H^3(K, \mu_2)$. Let $\calx$ be a projective regular model of
$K$. Let $C$ be an irreducible curve on $\calx$ and $v$ the
discrete valuation of $K$ given by $C$. Let $K_v$ be the
completion of $K$ at $v$. Then as in the proof of (4.4), we have
$u(K_v) = 8$. Thus by (4.2), there exist $g_v, h_v \in K_v^*$ such
that $g_v$ is a value of the quadratic form $q_3$ and $\phi = <\!
1, f\!><\! 1, g_v\!><\!1, h_v\!>$ over $K_v$. By weak
approximation, we can find $g \in K^*$ such that $g$ is a value of
$q_3$ over $K$ and $g = g_v$ modulo $K_v^{*^2}$ for all discrete
valuations $v$ of $K$ given by the irreducible curves $C$ in ${\rm
ram}_{\calx}(\zeta)$. Let $C$ be an irreducible curve on $\calx$
and $v$ the discrete valuation of $K$ given by $C$. By the choice
of $g$ it is clear that $\zeta = e_3(\phi) = (-f) \cdot (-g) \cdot
(-h_v)$ for all the discrete valuations $v$ of $K$ given by the
irreducible curves $C$ in the support of ram$_{\calx}(\zeta)$. If
$C$ is not in the support of ram$_{\calx}(\zeta)$, then as is the
proof of (4.4), we have $\zeta = (-f) \cdot (-g) \cdot (1)$ over
$K_v$. Let $\alpha = (-f) \cdot (-g) \in H^2(K, \mu_2)$. By (3.4),
there exists $h \in K^*$ such that $\zeta = \alpha \cup (-h) =
(-f) \cdot (-g) \cdot (-h)$. Since  $e_3 : I^3(K) \to H^3(K,
\mu_2)$ is an isomorphism, $ \phi = <\! 1, f\!><\! 1, g\!><\! 1,
h\!>$.

\paragraph*{Theorem 4.6} Let $K$ be a function field of a curve over
a $p$-adic field $k$. If $p \neq 2$, then $u(K) = 8$.

\paragraph*{Proof.} Let $K$ be a function field of a curve over a
$p$-adic field $k$. Assume that $p \neq 2$.    By a theorem of
Saltman ([S1], 3.4, cf. [S2]), every element in $H^2(K,\mu_2)$ is
a sum of at most 2 symbols. Since the cohomological dimension of
$K$ is 3, we also have $I^4(K) \simeq H^4(K, \mu_2) = 0$ ([AEJ]).
Now the theorem follows from (4.4), (4.5) and (4.3). \hfill
$\Box$.

\section*{References}

\begin{enumerate}

\item[{[A]}] Arason, J.K., \emph{Cohomologische Invarianten
quadratischer Formen, J. Algbera} {\bf 36} (1975), 448-491.

\item[{[AEJ]}] Arason, J.K., Elman, R. and Jacob, B., \emph{Fields
of cohomological \\ 2-dimension three, Math. Ann.} {\bf 274}
(1986),  649-657 .

\item[{[C]}] Colliot-Th\'el\`ene, J.-L., \emph{Birational
invariants, purity, and the Gresten conjecture, Proceedings of
Symposia in Pure Math.} {\bf 55}, Part 1, 1-64.

\item[{[K]}] Kato, K., \emph{A Hasse principle for two-dimensional
global fields, J. reine Angew. Math.} {\bf 366} (1986), 142-181.

\item[{[L]}] Lam, T.Y., \emph{Introduction to quadratic forms over
fields}, GSM 67, American Mathematical Society, 2004.

\item[{[Li1]}] Lipman, J., \emph{Introduction to resolution of
singularities, Proc. Symp. Pure Math.} {\bf 29} (1975), 187-230.

\item[{[Li2]}] Lipman, J., \emph{Desingularization  of
two-dimensional schemes, Ann. Math.} {\bf 107} (1978), 151-207.

\item[{[M]}] Merkurjev, A.S.,  \emph{On the norm residue symbol of
degree 2, Dokl. Akad. Nauk. SSSR} {\bf 261} (1981), 542-547.

\item[{[Mi]}] Milne, J.S.. \emph{\'Etale Cohomology}, Princeton
University Press, Princeton, New Jersey 1980.

\item[{[PS]}] Parimala, R. and Suresh, V., \emph{Isotropy  of
quadratic forms over function fields in one variable over $p$-adic
fields, Publ. de I.H.\'E.S.} {\bf 88} (1998), 129-150.

\item[{[S1]}] Saltman, D.J., \emph{Division Algebras over $p$-adic
curves,  J. Ramanujan Math. Soc.} {\bf 12} (1997), 25-47.

\item[{[S2]}] Saltman, D.J., \emph{Correction to Division algebras
over $p$-adic curves, J. Ramanujan Math. Soc.} {\bf 13} (1998),
125-130.

\item[{[S3]}] Saltman, D.J., Cyclic Algebras over $p$-adic curves
(to appear).

\item[{[Sc]}] Scharlau, W., \emph{Quadratic and Hermitian Forms},
Grundlehren der Math. Wiss., Vol. 270, Berlin, Heidelberg, New
York 1985.

\item[{[Sh]}] Shafarevich, I. R. \emph{Lectures on Minimal Models
and Birational Transformations of two Dimensional Schemes}, Tata
Institute of Fundamental Research, (1966).
\end{enumerate}

\noindent
Department of Mathematics and Computer Science \\
Emory University \\
400 Dowman Drive \\
Atlanta, Georgia 30322   \\
USA \\
E-mail: parimala@mathcs.emory.edu

\vskip 5mm

\noindent
Department of Mathematics and Statistics \\
University of Hyderabad \\
Gahcibowli \\
Hyderabad - 500046\\
India \\
E-mail: vssm@uohyd.ernet.in

\end{document}